\documentclass[runningheads]{llncs}
\usepackage[T1]{fontenc}
\usepackage{graphicx}
\usepackage{color}

\usepackage{hyperref}

\usepackage{lipsum}
\usepackage{amsfonts}
\usepackage{graphicx}
\usepackage{epstopdf}
\usepackage{bm}
\usepackage{dsfont}
\usepackage{enumitem}
\usepackage{amsmath}

\usepackage{algorithm}
\usepackage{algpseudocode}
\usepackage{amsmath}
\usepackage{xcolor}

\newcommand{\N}{{\mathbb{N}}}
\newcommand{\R}{{\mathbb{R}}}
\newcommand{\Cn}{{\mathbb{C}}}

\newcommand{\imag}{{\bm{\textrm{i}}}}

\newcommand{\rT}[1]{#1^{\textsf{T}}}

\newcommand{\Her}[1]{{{#1}^{\textsf{H}}}}
\newcommand{\Herb}[1]{{\Her{\left(#1\right)}}}

\newcommand{\si}[1]{{#1^{+}}}

\newcommand{\grad}[1][]{{\nabla_{#1}}}

\newcommand{\dd}[2][]{{\frac{\mathrm{d}^{#1}}{\mathrm{d}#2}}}

\newcommand{\ddt}{{\dd{t}}}

\newcommand{\Jtwo}[1]{{\mathbb{J}_{2#1}}}
\newcommand{\TJtwo}[1]{{\rT{\mathbb{J}}_{2#1}}}
\newcommand{\Jtn}{{\Jtwo{n}}}

\newcommand{\JtN}{{\Jtwo{N}}}
\newcommand{\TJtN}{{\TJtwo{N}}}

\newcommand{\I}[1]{{\fI_{#1}}}

\newcommand{\Z}[1]{{\fzero}_{#1}}

\newcommand{\tInit}{{t_{\mathrm{0}}}}
\newcommand{\tEnd}{{t_{\mathrm{end}}}}
\newcommand{\It}{I_t}

\newcommand{\fxInit}{{\fx_{\mathrm{0}}}}

\newcommand{\paramDomain}{\mathcal{P}}
\newcommand{\np}{{n_{\mathrm{\fmu}}}}

\newcommand{\reduce}[1]{{#1_{\mathrm{r}}}}
\newcommand{\fxr}{{\reduce{\fx}}}
\newcommand{\fxrInit}{{\fx}_{\mathrm{r,0}}}
\newcommand{\fxri}{{\fx}_{\mathrm{r},i}}

\newcommand{\ns}{{n_{\mathrm{s}}}}
\newcommand{\xs}{\fx^{\mathrm{s}}}

\newcommand{\fXs}{\fX_{\mathrm{s}}}
\newcommand{\fXsc}{\fXs^{\mathrm{c}}}

\newcommand{\Ham}{\mathcal{H}}

\newcommand{\Hamr}{\reduce{\Ham}}

\newcommand{\fe}{{\bm{e}}}

\newcommand{\fp}{{\bm{p}}}

\newcommand{\fv}{{\bm{v}}}

\newcommand{\fx}{{\bm{x}}}

\newcommand{\fA}{{\bm{A}}}
\newcommand{\fB}{{\bm{B}}}

\newcommand{\fD}{{\bm{D}}}

\newcommand{\fI}{{\bm{I}}}

\newcommand{\fK}{{\bm{K}}}

\newcommand{\fP}{{\bm{P}}}
\newcommand{\fQ}{{\bm{Q}}}

\newcommand{\fS}{{\bm{S}}}
\newcommand{\fT}{{\bm{T}}}
\newcommand{\fU}{{\bm{U}}}
\newcommand{\fV}{{\bm{V}}}
\newcommand{\fW}{{\bm{W}}}
\newcommand{\fX}{{\bm{X}}}
\newcommand{\fY}{{\bm{Y}}}

\newcommand{\fmu}{{\bm{\mu}}}

\newcommand{\fxi}{{\bm{\xi}}}

\newcommand{\fSigma}{{\bm{\Sigma}}}

\newcommand{\fOmega}{{\bm{\varOmega}}}

\newcommand{\povs}{p_\textrm{ovs}}
\newcommand{\qpow}{{q_\textrm{pow}}}

\newcommand{\fzero}{\ensuremath{\bm{0}}}

\newcommand{\calC}{\mathcal{C}}
\usepackage{pgfplots}
\usepgfplotslibrary{groupplots}
\usetikzlibrary{positioning}

\usepackage[capitalise,noabbrev]{cleveref}
\begin{document}
\title{Randomized Symplectic Model Order Reduction for Hamiltonian Systems}

\author{\underline{R. Herkert}\inst{1}\and
P. Buchfink\inst{1} \and
B. Haasdonk\inst{1}  \and 
J. Rettberg \inst{2} \and 
J. Fehr\inst{2}}
\authorrunning{R. Herkert et al.}
\institute{Institute of Applied Analysis and Numerical Simulation, University of Stuttgart, Pfaffenwaldring 57, 70569 Stuttgart, Germany \\
\email{\{robin.herkert,patrick.buchfink,haasdonk\}@mathematik.uni-stuttgart.de} 
\and Institute of Engineering and Computational Mechanics, University of Stuttgart, Pfaffenwaldring 9, 70569 Stuttgart, Germany \\
\email{\{johannes.rettberg,joerg.fehr\}@itm.uni-stuttgart.de}}

\maketitle           
\vspace{-7mm}
\begin{abstract}
Simulations of large scale dynamical systems in multi-query or real-time contexts require efficient surrogate modelling techniques, as e.g.\ achieved via  Model Order Reduction (MOR). Recently, symplectic methods like the complex singular value decomposition (cSVD) or the SVD-like decomposition have been developed for preserving Hamiltonian structure during MOR. 
In the current contribution, we show how symplectic structure preserving basis generation can be made more efficient with randomized matrix factorizations. We present a randomized complex SVD (rcSVD) algorithm and a randomized SVD-like (rSVD-like) decomposition. We demonstrate the efficiency of the approaches with numerical experiments on high dimensional systems. 

\end{abstract} \vspace{-7mm}

\section{Introduction} \vspace{-3mm}
Numerical simulation of large scale dynamical systems often leads to high computational costs. In a multi-query or real-time context, this requires efficient surrogate modelling techniques such as model order reduction (MOR). Additionally, models appear as Hamiltonian systems, which, for example, describe conservative dynamics and non-dissipative phenomena in classical mechanics or transport problems.
The structure of a Hamiltonian system ensures conservation of energy and, under mild assumptions, stability properties. Classical MOR like the Proper Orthogonal Decomposition (POD)  \cite{Volkwein2013} fails to preserve this Hamiltonian structure which, in general, violates the conservation of energy and may yield unstable reduced models. Thus, one current trend in MOR is the structure-preserving MOR of Hamiltonian systems \cite{Buchfink2019,MaboudiAfkham2017}. Recently, structure-preserving methods like the complex singular value decomposition (cSVD) \cite{Peng2016} or the SVD-like decomposition \cite{Buchfink2019,Xu2003} have been developed. Both of these methods belong to the class of data-based/snapshot-based MOR, which has the advantage that it can (i) also be used in black-box learning, (ii) be applied to closed-source nonlinear models where the underlying structure of the model is not clear. 
\\
With projection-based MOR (e.g.\ POD, cSVD and SVD-like), the solution is approximated in a low-dimensional subspace. For computing a basis of such a subspace, low-rank matrix approximations, like the truncated singular value decomposition (SVD), are computed for the snapshot matrix. In many cases, a randomized approach for a low-rank matrix approximation is faster and/or more robust than its classical version \cite{Halko2011}. How randomization can be applied to MOR is currently intensively studied: Randomized versions of classical MOR basis generation algorithms have been recently applied \cite{Alla2019,Bach2019,Buhr2018} or the reduced order model itself is approximated by its random sketch \cite{Balabanov2019a,Balabanov2019}. It has been shown that randomization can improve not only efficiency but also numerical stability. None of these approaches, in general, preserves a Hamiltonian structure of a system. 
 \\
This work is focussed on how structure-preserving symplectic basis generation and efficient randomized basis generation can be combined. We present two randomized, symplectic basis generation schemes: a randomized complex SVD (rcSVD) algorithm and a randomized SVD-like (rSVD-like) decomposition. For the rcSVD, we leverage randomization for complex matrices. The rSVD-like algorithm is obtained by a randomized version of the Schur decomposition. We demonstrate the efficiency of the approaches by numerical testing on high-dimensional systems such as obtained from spatial discretization of the wave equation. This work is structured as follows: In \Cref{sec:essentials} essential background on MOR, Hamiltonian systems and randomized matrix factorizations is given. In \Cref{sec:main}, the new randomized, symplectic methods are presented. \Cref{sec:numerics} is focussed on numerical experiments and comparisons with non-randomized structure-preserving methods.
The work is concluded in \Cref{sec:conclusion}. 

\vspace{-4mm}
\section{Essentials}\vspace{-3mm}
\label{sec:essentials}
\textbf{Symplectic Model Order Reduction of Hamiltonian Systems} \vspace{1mm} \\
In this section, a brief summary on MOR of Hamiltonian systems is given. For a more detailed introduction, we refer to \cite{Benner2017} (MOR),  \cite{daSilva2008} (symplectic geometry and Hamiltonian systems) and \cite{Peng2016} (symplectic MOR of Hamiltonian systems). \\ 
Given a parametric \emph{Hamiltonian (function)} $\Ham(\cdot; \fmu) \in \calC^1(\R^{2N})$ and a parametric initial value $\fxInit(\fmu)\in \R^{2N}$, a parametric, canonical \emph{Hamiltonian system} reads: For a fixed (but arbitrary) parameter vector
$\fmu \in \paramDomain \subset \R^{\np}$ and time interval $I_t := [\tInit, \tEnd]$\footnote{Note that $I_t$ may even be parameter dependent, as we will use in the experiments. For notational simplicty we keep a fixed time interval in this section.}, find the full solution $\fx(\cdot; \fmu) \in \calC^1(\It, \R^{2N})$
with
\begin{equation}\label{eq:ham_sys}
  \ddt \fx(t; \fmu)
= \JtN \grad[\fx] \Ham(\fx(t; \fmu); \fmu), \quad   \fx(\tInit; \fmu) = \fxInit(\fmu)
\end{equation}
\vspace{-3mm}
where  $$
  \JtN := \begin{bmatrix}
    \Z{N} & \I{N}\\
    -\I{N} & \Z{N}
  \end{bmatrix} \in \R^{2N \times 2N}$$ 
denotes the canonical Poisson matrix and $\I{N}, \Z{N}\in \R^{N\times N}$ denote the identity and zero matrices. 
The most important property of Hamiltonian systems is
that the solution conserves the Hamiltonian over time, i.e.
$\ddt \Ham(\fx(t; \fmu); \fmu) = 0$ for all $t \in \It$.
Given a reduced order basis (ROB) matrix $\fV\in\R^{2 N \times 2n}$ and projection matrix $\fW\in\R^{2 N \times 2n}$, with $\rT\fW\fV = \I{2n}$ and $n \ll N,$
the projection-based reduction of the Hamiltonian system leads to a reduced system that reads: 
For a parameter vector $\fmu \in \paramDomain \subset \R^{\np}$,
find the \emph{reduced solution} $\fxr(\cdot; \fmu) \in \calC^1(\It, \R^{2n})$ 
with \vspace{-3mm}
\begin{equation}\label{eq:ROM}
  \ddt \fxr(t; \fmu) = \rT\fW\JtN \grad[\fx]\Ham(\fV\fxr(t; \fmu); \fmu), \quad  \fxr(\tInit; \fmu) = \fxrInit(\fmu)
\end{equation}
with the \emph{reduced initial value}
$\fxrInit(\fmu) := \rT\fW \fxInit(\fmu)$. This system is not necessarily Hamiltonian and thus the conservation of the Hamiltonian over time can not be ensured. 
To preserve the Hamiltonian structure, symplectic MOR can be used \cite{Buchfink2019,Peng2016}. Here, (i) the ROB matrix $\fV$ is required to be a \emph{symplectic matrix}
which means, that for $\fV \in \R^{2N \times 2n}$ with $n \leq N$
\begin{align*}
  \rT\fV \JtN \fV = \Jtn,
\end{align*}
and (ii) the projection matrix $\fW$ is set to be the transpose of the so-called \emph{symplectic inverse}
$\si{\fV}$ of the ROB matrix $\fV$, i.e.
\begin{align*}
  \rT\fW = \si{\fV} := \Jtn \rT\fV \TJtN.
\end{align*}

With this choice, the reduced model \eqref{eq:ROM} is a low-dimensional
Hamiltonian system with the \emph{reduced Hamiltonian} $\Hamr(\fxr(t;\fmu);\fmu) := \Ham(\fV \fxr(t;\fmu);\fmu)$, which is defined as the Hamiltonian of the ROB matrix times the reduced coordinates.\\ \\
\textbf{Randomized Matrix Factorizations} \vspace{1mm} \\
We continue with a brief summary on randomized matrix factorizations. For a more detailed presentation, we refer to \cite{Halko2011}.
The computation of a randomized factorization of a matrix $\fB \in \R^{m \times l}$ is subdivided into two stages. First, for $k \leq m$ a matrix $\fQ \in \R^{m \times k}$ with orthonormal columns is computed that approximates $\fB \approx \fQ \rT\fQ\fB$. This task can be efficiently executed with random sampling methods. 
Then, a matrix decomposition (e.g.\ SVD, QR) of $\rT\fQ \fB \in \R^{k \times l}$ is computed and multiplied by $\fQ$. The first factor of a SVD/QR-decomposition is a matrix with orthonormal columns and this property is not changed by a multiplication with $\fQ$. For the SVD, we get $\rT\fQ \fB = \widetilde\fU\fSigma\rT\fV$ and by setting $\fU:= \fQ \widetilde \fU$ we get the approximate SVD $ \fB \approx \fU\fSigma\rT\fV.$ Instead of the target rank $k$, it is known to be advantageous to introduce an oversampling parameter $\povs$ and aim for  $k+\povs$ columns for $\fQ$, then for $\fU$ only the first $k$ columns are used. 
The matrix $\fQ$ is computed via \Cref{randomsampling}. Note, that a computational advantage over direct factorization of $\fB$ will be achievable if $k\ll l$.
\begin{algorithm}[H]
\caption{Random Sampling Algorithm}
\label{randomsampling}
{\bf Input}: $\fB\in \R^{m \times l}$, target rank $k \in \N$, oversampling parameter
$\povs \in \N_0$\\
{\bf Output}: Matrix with orthonormal columns $\fQ\in \R^{m \times (k+\povs)}$
\begin{algorithmic}[1]
\State Draw a random Gaussian test matrix $\fOmega \in \R^{l \times (k + \povs)}$.
\State Compute the matrix product $\fY = \fB\fOmega$, a so-called \emph{random sketch} of $\fB$.
\State Construct a matrix $\fQ\in \R^{m \times (k+\povs)}$ whose columns form an orthonormal basis for the range of $\fY$\footnotemark.
\end{algorithmic}
\end{algorithm}
\footnotetext{In the rare case of linear dependencies in $\fY$, the number of columns in $\fQ$ is reduced.}
Step 3 of this algorithm can for example be performed using a QR-decom-position or a SVD. Using a random matrix $\fOmega$, with a special factorization such as SRFT \cite[Section~4.6]{Halko2011}, randomized schemes can produce an approximate SVD using only
$\mathcal{O}(ml\log(k+\povs) + (m + l)(k+\povs)^2)$ flops, because the multiplication $\fB\fOmega$ can be computed in $\mathcal{O}(ml\log(k+\povs))$. 
In contrast, the cost of a classical approach is typically $\mathcal{O}(mlk)$ flops.
For the projection error $||\fB-\fQ\rT\fQ\fB||_2$, a probabalistic error bound can be proven: With probability at least $1-3\povs^{-\povs}$ the bound $||\fB-\fQ\rT\fQ\fB||_2\leq (1 + \sqrt{9 k + \povs}$ min$(m, l))\sigma_{k+1}$
holds under mild assumptions on $\povs$, with $\sigma_{k+1}$ denoting the $(k+1)$th singular value of $\fB$ \cite[Section 10.3]{Halko2011}. 
This error (bound) can be further improved by a power iteration, which means that for $\qpow \in \N_0$, the random sketch is computed as $\fB(\rT\fB\fB)^{\qpow}\fOmega$. This is in particular useful for matrices whose singular values decay slowly. Also, randomized a posteriori error estimation is possible, even at almost no additional computational cost (see \cite[Sections~4.3 and 4.4]{Halko2011}).  

\vspace{-4mm}
\section{Randomized Symplectic Model Order Reduction}\vspace{-3mm}
\label{sec:main}
In this section, two new randomized symplectic methods are presented.
In the following $\imag \in \Cn$ denotes the imaginary unit and $\Herb{\cdot}$ the complex transpose.
Moreover, we use MATLAB-style notation for matrix indexing and stacking.\\ \\
\textbf{Randomized Complex SVD}\vspace{1mm} \\
Consider the snapshot matrix $\fXs := (\xs_i)_{i=1}^{\ns} \in \R^{2N \times \ns}$, where $\fXs$ is split into $\fXs = [\fQ_s; \fP_s]$, with $\fQ_s, \fP_s \in \R^{N \times \ns}$.
The main idea of the cSVD algorithm is to form a complex snapshot matrix $\fXsc := \fQ_s + \imag \fP_s \in \Cn^{N \times \ns}$ and compute a truncated SVD of this complex matrix
$\fXsc \approx \fU_\textrm{C} \fSigma_\textrm{C}\rT\fV_\textrm{C}.$ The matrix $\fU_\textrm{C} \in \Cn^{N \times k}$ is then split into real and imaginary part $\fU_\textrm{C} = \fV_Q + \imag \fV_P$ and mapped to $$\fV := \mathcal{A}(\fU_\textrm{C}) := 
\begin{pmatrix}
&\fV_Q \  &-\fV_P \\
&\fV_P \  &\fV_Q
\end{pmatrix}.$$
This mapping $\mathcal {A}$ from the complex Stiefel manifold $\fV_k(\Cn^{n})$ to $\R^{2N\times2k}$ maps a complex matrix with orthonormal columns to a real symplectic matrix (see \cite{Peng2016}). The symplectic matrix $\fV$ and its symplectic inverse $ \si{\fV}$ are then used for MOR. 
Instead of using a truncated SVD, we apply randomization in order to compute a rank-$k$ approximation of $\fXsc.$ 
The procedure is summarized as \cref{rcSVDalg}.
\begin{algorithm}
\caption{Randomized Complex SVD (rcSVD)}\label{rcSVDalg}
{\bf Input}: Snapshot matrix $\fXs\in \R^{2N\times n_s}$ ,
target rank $2k \in \N$ of the ROB,
oversampling parameter $\povs \in \N_0$,
power iteration number $\qpow \in \N_0$\\
{\bf Output}: Symplectic ROB matrix $\fV_\textrm{rcSVD} \in \R^{2N\times 2k}$
\begin{algorithmic}[1]
\State $\fXsc = \fXs(1 : N, :) +  \imag  \fXs((N + 1) : (2N), :)$\Comment{construct complex snapshot matrix}
\State $\fOmega = \mathtt{SRFT}(N, k+\povs)$\Comment{draw a random $N \times (k + \povs)$ test matrix}
\State $[\fU_\textrm{C},\fSigma_\textrm{C}, \fV_\textrm{C}] = \mathtt{SVD}(\fXsc(\Herb\fXsc \fXsc)^{\qpow}\fOmega$, $k$) \Comment{basis for approximation of $\fY_c$}
\State$\fV_\textrm{Q} =$ Re($\fU_\textrm{C}$) , $\fV_\textrm{P}$ = Im($\fU_\textrm{C}$)  \Comment{split in real and imaginary part}
\State$\fV_{\textrm{rcSVD}} = [\fV_\textrm{Q}, -\fV_\textrm{P}; \fV_\textrm{P}, \fV_\textrm{Q}]$  \Comment{map to symplectic matrix}
\end{algorithmic}
\end{algorithm} \\ \\
\textbf{Randomized SVD-like}\vspace{1mm} \\
In \cite{Xu2003} it is shown that each real $2N \times \ns$ matrix can be decomposed as
$\fXs = \fS\fD\rT\fP$, with $\fS\in \R^{2N \times 2N}$ symplectic, $\fP\in  \R^{\ns\times \ns}$ orthogonal,
\[ \begin{array}{r@{\,}l}
    & 
    \begin{matrix}
      \mspace{15mu}&\overbrace{\rule{0.2cm}{0pt}}^{p} & \overbrace{\rule{0.2cm}{0pt}}^{q} &\overbrace{\rule{0.2cm}{0pt}}^{N-p-q} &\overbrace{\rule{0.2cm}{0pt}}^{p} & \overbrace{\rule{0.2cm}{0pt}}^{N-p}
    \end{matrix}
    \\
    \rT\fD = & \begin{pmatrix}
&\fSigma \hspace{0.5cm}  &\Z{} \hspace{0.5cm}  &\Z{} \hspace{0.5cm}  &\Z{} \hspace{0.5cm}  &\Z{} \\
&\Z{} \hspace{0.5cm} &\I{q} \hspace{0.5cm} &\Z{} \hspace{0.5cm}&\Z{} \hspace{0.5cm} &\Z{} \\	
&\Z{} \hspace{0.5cm} &\Z{} \hspace{0.5cm} &\Z{} \hspace{0.5cm} &\fSigma \hspace{0.5cm}  &\Z{}\\
&\Z{} \hspace{0.5cm} &\Z{} \hspace{0.5cm} &\Z{} \hspace{0.5cm} &\Z{}\hspace{0.5cm}	 &\Z{}  \\	
\end{pmatrix}	
  \end{array}
 \in \R^{\ns\times 2N},\quad
\fSigma = 
\begin{pmatrix} 
&\sigma_1 \ &0 \ &0 \\
&0 \ & \ddots \ &0  \\	
&0 \ &0  \ &\sigma_p \\
\end{pmatrix}
\in \R^{p \times p}\]
with $\sigma_i >0$, $i= 1,\dots,p$. Note that due to symplecticity of $\fS$
$$\fK : =\rT\fXs\JtN\fXs = \fP\rT\fD\rT\fS\JtN\fS\fD\rT\fP =\fP\rT\fD\JtN\fD\rT\fP,$$ where $\rT\fD\JtN\fD \in \R^{\ns \times \ns}$ is a matrix with $\left\{ \pm\sigma_i^2 \right\}_{i=1}^p$ on the $(p+q)$th superdiagonal and subdiagonal and zeros everywhere else.  Thus, the factorization $\fP\rT\fD\JtN\fD\rT\fP$ is a permutation of the real Schur decomposition $\fK =  \fU \fT  \rT \fU$, where $\left\{+\sigma_i^2 \right\}_{i=1}^p$ are on the first superdiagonal and $\left\{ -\sigma_i^2 \right\}_{i=1}^p$ on the first subdiagonal of $\fT$ and $\fU$ is orthogonal. Instead of performing a standard Schur decomposition, we randomize the Schur decomposition of $\fK$ and obtain a randomized SVD-like decomposition in this way. The procedure is presented as \cref{rSVDlikealg}. For the function $\mathtt{computeQ}(\povs, \qpow, \fXs)$ in Step 1, one of the following methods has to be inserted, where $\mathtt{randn}$ generates a random Gaussian matrix:
\begin{enumerate}
\label{enumcompQ}
  \item $\mathtt{computeQfromK}(\povs, \qpow, \fXs)$: \\
$\fK = \rT\fXs\JtN\fXs, \fOmega = \mathtt{randn}(\ns, k+\povs)$, $\fQ =\mathtt{orth}(\fK^{2\qpow+1} \fOmega$)
  \item $\mathtt{computeQfromXs}(\povs, \qpow, \fXs)$: \\
$\fOmega = \mathtt{randn}(2N, k+\povs)$, $\fQ =\mathtt{orth}(\rT\fXs (\fXs\rT\fXs)^\qpow\fOmega$)
  \item $\mathtt{computeQfromKXs}(\povs, \qpow, \fXs)$:  \\
$\fOmega_K = \mathtt{randn}(\ns,\lceil\frac{k+\povs}{2}\rceil)$, $\fOmega_X = \mathtt{randn}(2N, \lfloor\frac{k+\povs}{2}\rfloor)$, \vspace{1mm}
  \item[] $\fK = \rT\fXs\JtN\fXs, \fQ =\mathtt{orth}( [\fK^{2\qpow+1 }\fOmega_K,\rT\fXs  (\fXs\rT\fXs)^\qpow\fOmega_X]$).
\end{enumerate}
The resulting variants of the rSVD-like algorithm will be named accordingly: rSVD-like\texttt{*}, where \texttt{*} is to be replaced by \texttt{K}, \texttt{Xs} or \texttt{KXs} depending on the computational variant for $\fQ$.

\vspace{-4mm}
\begin{algorithm}[H]
\caption{Randomized SVD-like (rSVD-like)}\label{rSVDlikealg}
{\bf Input}:
Snapshot matrix $\fXs\in \R^{2N\times n_s}$,
target rank $2k \in \N$ of the ROB,
oversampling parameter $\povs \in \N_0$,
power iteration number $\qpow \in \N_0$\\
{\bf Output}: Symplectic ROB matrix $\fV_\textrm{rSVD-like} \in \R^{2N\times 2k}$
\begin{algorithmic}[1]
\State $\fQ = \mathtt{computeQ}(\povs, \qpow, \fXs)$ \Comment{see end \Cref{enumcompQ}}
\State $\fK = \rT\fXs\JtN\fXs$
\State $[\fU, \fT] = \mathtt{realSchur}(\rT\fQ\fK\fQ$) \Comment{compute real Schur decomposition}
\State $p =$ rank$(\fT)/2$
\State $\fSigma = \mathtt{diag}(\sqrt{\fT_{1,2}}, \sqrt{ \fT_{3,4}},...,  \sqrt{\fT_{2p-1,2p}}$) \Comment{extract real $\sigma_i, i = 1,..,p$}
\State $\fP :=[\fp_1, ..., \fp_{\ns}] := \fU\cdot[\I{\ns}$(:,1:2:2$p$-1), $\I{\ns}$(:,2:2:2$p$), $\I{\ns}$(:,2$p$+1:$\ns$)]
\State $\fV_\textrm{rSVD-like} = \fXs\ [\fP(:,1:k)\fSigma(1:k,1:k)^{-1},\fP(:,p+1:p+k)\fSigma(1:k,1:k)^{-1} ]$
\end{algorithmic}
\end{algorithm}
\vspace{-4mm}
\section{Numerical experiments}
\begin{figure}[h]
  \centering
\pgfplotsset{
  every axis plot/.append style={line width=1pt,},
  every axis/.append style={
    ymajorgrids,
    xmajorgrids,
    grid style={dashed, lightgray,semithick},
   	axis line style = semithick,
   	every tick/.style={semithick,},
	  yticklabels={,,},
  },
}
\pgfplotsset{
  /pgfplots/group/every plot/.append style = {
    height = .3\textwidth, width = .34\textwidth,
    ylabel near ticks,
    legend style={
      draw=none,
      fill=none,
      legend cell align=left,
      column sep=.5em,
    },
  },
}
\pgfplotsset{%
my legend/.style={legend image code/.code={%
\node[##1,anchor=west] at (0cm,0cm){\pgfuseplotmark{+}};
\node[##1] at (0.4cm,0cm){\pgfuseplotmark{*}};
}},%
relerr/.style={
  xmin=10, xmax=1000,
  ymin=4e-7, ymax=1e0,
  ytick={1e-6, 1e-4, 1e-2, 1e0},
},
runtime/.style={
  xmin=10, xmax=1000,
  ymin=5e-1, ymax=1000,
  ytick={1e0, 1e1, 1e2, 1e3},
},
cSVDwSVDs/.style={color0,dashed},
cSVDev/.style={color1,dashed},
rSVDlikep0q0/.style={color5,dashed},
rSVDlikep10q0/.style={color3,dashed},
rSVDlikep10q1/.style={color4,dashed},
cSVDfull/.style={color2,dashed},
row0/.style={xticklabels={,,},},
row1/.style={xlabel={basis size}, xlabel near ticks},
}

\definecolor{color0}{rgb}{0.12156862745098,0.466666666666667,0.705882352941177}
\definecolor{color1}{rgb}{0.83921568627451,0.152941176470588,0.156862745098039}
\definecolor{color2}{rgb}{1,0.498039215686275,0.0549019607843137}
\definecolor{color3}{rgb}{0.549019607843137,0.337254901960784,0.294117647058824}
\definecolor{color4}{rgb}{0.890196078431372,0.466666666666667,0.76078431372549}
\definecolor{color5}{rgb}{0.580392156862745,0.403921568627451,0.741176470588235}
\definecolor{color6}{rgb}{0.172549019607843,0.627450980392157,0.172549019607843}

\begin{tikzpicture}
\begin{groupplot}[
  group style={
    group name=err_plots,
    group size=2 by 1,
    vertical sep=1cm,
    horizontal sep=1cm
}]
\nextgroupplot[
    runtime,
    row1,
    ymode=log,
    xmode=log,
    legend to name=leg0,
    legend columns = 1,
    title={basis gen. time},
    yticklabels={1e0, 1e1, 1e2, 1e3},
]
  \addplot[cSVDwSVDs,dash phase=2,forget plot] table [
    header=true,
    col sep=comma,
    x=rbsize,
    y=runtime,
  ]{pics/data/plot_cSVDwithsvds.dat};
  \addlegendimage{cSVDwSVDs} 
  \addlegendentry{cSVD with \texttt{SVDs}}

  \addplot[cSVDev] table [
    header=true,
    col sep=comma,
    x=rbsize,
    y=runtime,
  ]{pics/data/plot_cSVDev.dat};
  \addlegendentry{cSVD with \texttt{eig}}

  \addplot[rSVDlikep0q0,dash phase=3,forget plot] table [
    header=true,
    col sep=comma,
    x=rbsize,
    y=runtime,
  ]{pics/data/plot_rcSVD_povs0_qpow0.dat};
  \addlegendimage{rSVDlikep0q0} 
  \addlegendentry{rcSVD, $p_{\text{ovs}} = 0$, $q_{\text{pow}} = 0$}

  \addplot[rSVDlikep10q0,dash phase=2,forget plot] table [
    header=true,
    col sep=comma,
    x=rbsize,
    y=runtime,
  ]{pics/data/plot_rcSVD_povs10_qpow0.dat};
  \addlegendimage{rSVDlikep10q0} 
  \addlegendentry{rcSVD, $p_{\text{ovs}} = 10$, $q_{\text{pow}} = 0$}

  \addplot[rSVDlikep10q1] table [
    header=true,
    col sep=comma,
    x=rbsize,
    y=runtime,
  ]{pics/data/plot_rcSVD_povs10_qpow1.dat};
  \addlegendentry{rcSVD, $p_{\text{ovs}} = 10$, $q_{\text{pow}} = 1$}

  \addplot[cSVDfull] table [
    header=true,
    col sep=comma,
    x=rbsize,
    y=runtime,
  ]{pics/data/plot_cSVDfull.dat};
  \addlegendentry{cSVD full}
  \nextgroupplot[
      relerr,
      row1,
      ymode=log,
      xmode=log,
      title={rel. reduction error},
      yticklabels={1e-6,1e-4,1e-2,1e0},
  ]
  \addplot[cSVDwSVDs,dash phase=2] table [
    header=true,
    col sep=comma,
    x=rbsize,
    y=relerr,
  ]{pics/data/plot_cSVDwithsvds.dat};

  \addplot[cSVDev] table [
    header=true,
    col sep=comma,
    x=rbsize,
    y=relerr,
  ]{pics/data/plot_cSVDev.dat};

  \addplot[rSVDlikep0q0,dash phase=3] table [
    header=true,
    col sep=comma,
    x=rbsize,
    y=relerr,
  ]{pics/data/plot_rcSVD_povs0_qpow0.dat};

  \addplot[rSVDlikep10q0,dash phase=2] table [
    header=true,
    col sep=comma,
    x=rbsize,
    y=relerr,
  ]{pics/data/plot_rcSVD_povs10_qpow0.dat};

  \addplot[rSVDlikep10q1] table [
    header=true,
    col sep=comma,
    x=rbsize,
    y=relerr,
  ]{pics/data/plot_rcSVD_povs10_qpow1.dat};

  \addplot[cSVDfull] table [
    header=true,
    col sep=comma,
    x=rbsize,
    y=relerr,
  ]{pics/data/plot_cSVDfull.dat};
\end{groupplot}

\node[
    right = .35cm of err_plots c2r1.east,
    align=left,
    draw=black,
    anchor=west,
  ] {%
  \scriptsize
  \ref*{leg0}
  };
\end{tikzpicture}%
  \caption{Basis generation times and relative reduction errors (\ref{eqnrelerr}), rcSVD}%
  \label{fig:methcomp2}%
\end{figure}
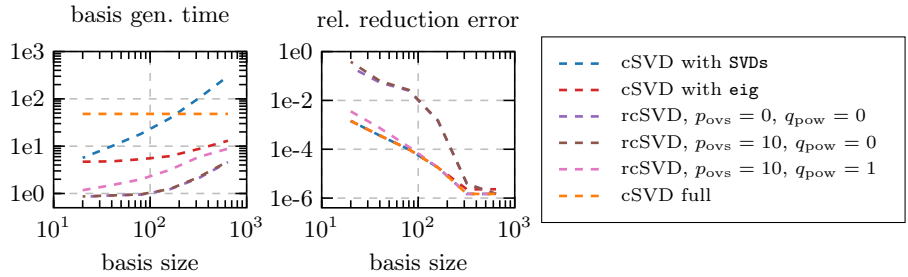

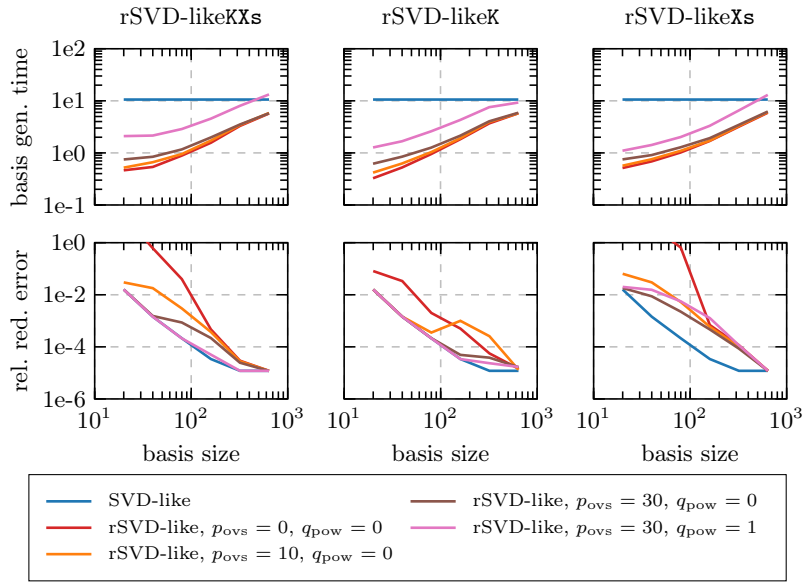
\begin{figure}[H]
  \centering
\pgfplotsset{
  every axis plot/.append style={line width=1pt,},
  every axis/.append style={
    ymajorgrids,
    xmajorgrids,
    grid style={dashed, lightgray,semithick},
   	axis line style = semithick,
   	every tick/.style={semithick,},
	  yticklabels={,,},
  },
}
\pgfplotsset{
  /pgfplots/group/every plot/.append style = {
    height = .3\textwidth, width = .34\textwidth,
    ylabel near ticks,
    legend style={
      draw=none,
      fill=none,
      legend cell align=left,
      column sep=.5em,
    },
  },
}
\pgfplotsset{%
my legend/.style={legend image code/.code={%
\node[##1,anchor=west] at (0cm,0cm){\pgfuseplotmark{+}};
\node[##1] at (0.4cm,0cm){\pgfuseplotmark{*}};
}},%
relerr/.style={
  xmin=10, xmax=1000,
  ymin=0.000001, ymax=1e0,
  ytick={1e-6, 1e-4, 1e-2, 1e0},
},
runtime/.style={
  xmin=10, xmax=1000,
  ymin=1e-1, ymax=100,
  ytick={1e-1, 1e0, 1e1, 1e2},
},
SVDlike/.style={color0,},
rSVDlikep0q0/.style={color1,},
rSVDlikep10q0/.style={color2,},
rSVDlikep30q0/.style={color3},
rSVDlikep30q1/.style={color4},
pod_proj/.style={color=black, dashed, mark=none},
row0/.style={xticklabels={,,},},
row1/.style={xlabel={basis size}, xlabel near ticks},
}

\definecolor{color0}{rgb}{0.12156862745098,0.466666666666667,0.705882352941177}
\definecolor{color1}{rgb}{0.83921568627451,0.152941176470588,0.156862745098039}
\definecolor{color2}{rgb}{1,0.498039215686275,0.0549019607843137}
\definecolor{color3}{rgb}{0.549019607843137,0.337254901960784,0.294117647058824}
\definecolor{color4}{rgb}{0.890196078431372,0.466666666666667,0.76078431372549}
\definecolor{color5}{rgb}{0.580392156862745,0.403921568627451,0.741176470588235}
\definecolor{color6}{rgb}{0.172549019607843,0.627450980392157,0.172549019607843}

\begin{tikzpicture}
\begin{groupplot}[
  group style={
    group name=err_plots,
    group size=3 by 2,
    vertical sep=.5cm,
    horizontal sep=.75cm
}]
\nextgroupplot[
    runtime,
    row0,
    ymode=log,
    xmode=log,
    legend to name=leg0,
    legend columns = 3,
    transpose legend,
    title={rSVD-like\texttt{KXs}},
    yticklabels={1e-1, 1e0, 1e1, 1e2},
]
  \addplot[SVDlike] table [
    header=true,
    col sep=comma,
    x=rbsize,
    y=runtime,
  ]{pics/data/plot_SVD-like.dat};
  \addlegendentry{SVD-like}

  \addplot[rSVDlikep0q0] table [
    header=true,
    col sep=comma,
    x=rbsize,
    y=runtime,
  ]{pics/data/plot_rSVD-likeKXS_povs0_qpow0.dat};
  \addlegendentry{rSVD-like, $p_{\text{ovs}} = 0$, $q_{\text{pow}} = 0$}

  \addplot[rSVDlikep10q0] table [
    header=true,
    col sep=comma,
    x=rbsize,
    y=runtime,
  ]{pics/data/plot_rSVD-likeKXS_povs10_qpow0.dat};
  \addlegendentry{rSVD-like, $p_{\text{ovs}} = 10$, $q_{\text{pow}} = 0$}

  \addplot[rSVDlikep30q0] table [
    header=true,
    col sep=comma,
    x=rbsize,
    y=runtime,
  ]{pics/data/plot_rSVD-likeKXS_povs30_qpow0.dat};
  \addlegendentry{rSVD-like, $p_{\text{ovs}} = 30$, $q_{\text{pow}} = 0$}

  \addplot[rSVDlikep30q1] table [
    header=true,
    col sep=comma,
    x=rbsize,
    y=runtime,
  ]{pics/data/plot_rSVD-likeKXS_povs30_qpow1.dat};
  \addlegendentry{rSVD-like, $p_{\text{ovs}} = 30$, $q_{\text{pow}} = 1$}

  \nextgroupplot[
    runtime,
    row0,
    ymode=log,
    xmode=log,
    title={rSVD-like\texttt{K}},
]
  \addplot[SVDlike] table [
    header=true,
    col sep=comma,
    x=rbsize,
    y=runtime,
  ]{pics/data/plot_SVD-like.dat};

  \addplot[rSVDlikep0q0] table [
    header=true,
    col sep=comma,
    x=rbsize,
    y=runtime,
  ]{pics/data/plot_rSVD-likeK_povs0_qpow0.dat};

  \addplot[rSVDlikep10q0] table [
    header=true,
    col sep=comma,
    x=rbsize,
    y=runtime,
  ]{pics/data/plot_rSVD-likeK_povs10_qpow0.dat};

  \addplot[rSVDlikep30q0] table [
    header=true,
    col sep=comma,
    x=rbsize,
    y=runtime,
  ]{pics/data/plot_rSVD-likeK_povs30_qpow0.dat};

  \addplot[rSVDlikep30q1] table [
    header=true,
    col sep=comma,
    x=rbsize,
    y=runtime,
  ]{pics/data/plot_rSVD-likeK_povs30_qpow1.dat};

  \nextgroupplot[
    runtime,
    row0,
    ymode=log,
    xmode=log,
    title={rSVD-like\texttt{Xs}},
]
  \addplot[SVDlike] table [
    header=true,
    col sep=comma,
    x=rbsize,
    y=runtime,
  ]{pics/data/plot_SVD-like.dat};

  \addplot[rSVDlikep0q0] table [
    header=true,
    col sep=comma,
    x=rbsize,
    y=runtime,
  ]{pics/data/plot_rSVD-likeXs_povs0_qpow0.dat};

  \addplot[rSVDlikep10q0] table [
    header=true,
    col sep=comma,
    x=rbsize,
    y=runtime,
  ]{pics/data/plot_rSVD-likeXs_povs10_qpow0.dat};

  \addplot[rSVDlikep30q0] table [
    header=true,
    col sep=comma,
    x=rbsize,
    y=runtime,
  ]{pics/data/plot_rSVD-likeXs_povs30_qpow0.dat};

  \addplot[rSVDlikep30q1] table [
    header=true,
    col sep=comma,
    x=rbsize,
    y=runtime,
  ]{pics/data/plot_rSVD-likeXs_povs30_qpow1.dat};
  \nextgroupplot[
      relerr,
      row1,
      ymode=log,
      xmode=log,
      yticklabels={1e-6,1e-4,1e-2,1e0},
  ]
    \addplot[SVDlike] table [
      header=true,
      col sep=comma,
      x=rbsize,
      y=relerr,
    ]{pics/data/plot_SVD-like.dat};

    \addplot[rSVDlikep0q0] table [
      header=true,
      col sep=comma,
      x=rbsize,
      y=relerr,
    ]{pics/data/plot_rSVD-likeKXS_povs0_qpow0.dat};

    \addplot[rSVDlikep10q0] table [
      header=true,
      col sep=comma,
      x=rbsize,
      y=relerr,
    ]{pics/data/plot_rSVD-likeKXS_povs10_qpow0.dat};

    \addplot[rSVDlikep30q0] table [
      header=true,
      col sep=comma,
      x=rbsize,
      y=relerr,
    ]{pics/data/plot_rSVD-likeKXS_povs30_qpow0.dat};

    \addplot[rSVDlikep30q1] table [
      header=true,
      col sep=comma,
      x=rbsize,
      y=relerr,
    ]{pics/data/plot_rSVD-likeKXS_povs30_qpow1.dat};

    \nextgroupplot[
      relerr,
      row1,
      ymode=log,
      xmode=log,
  ]
    \addplot[SVDlike] table [
      header=true,
      col sep=comma,
      x=rbsize,
      y=relerr,
    ]{pics/data/plot_SVD-like.dat};

    \addplot[rSVDlikep0q0] table [
      header=true,
      col sep=comma,
      x=rbsize,
      y=relerr,
    ]{pics/data/plot_rSVD-likeK_povs0_qpow0.dat};

    \addplot[rSVDlikep10q0] table [
      header=true,
      col sep=comma,
      x=rbsize,
      y=relerr,
    ]{pics/data/plot_rSVD-likeK_povs10_qpow0.dat};

    \addplot[rSVDlikep30q0] table [
      header=true,
      col sep=comma,
      x=rbsize,
      y=relerr,
    ]{pics/data/plot_rSVD-likeK_povs30_qpow0.dat};

    \addplot[rSVDlikep30q1] table [
      header=true,
      col sep=comma,
      x=rbsize,
      y=relerr,
    ]{pics/data/plot_rSVD-likeK_povs30_qpow1.dat};

    \nextgroupplot[
      relerr,
      row1,
      ymode=log,
      xmode=log,
  ]
    \addplot[SVDlike] table [
      header=true,
      col sep=comma,
      x=rbsize,
      y=relerr,
    ]{pics/data/plot_SVD-like.dat};

    \addplot[rSVDlikep0q0] table [
      header=true,
      col sep=comma,
      x=rbsize,
      y=relerr,
    ]{pics/data/plot_rSVD-likeXs_povs0_qpow0.dat};

    \addplot[rSVDlikep10q0] table [
      header=true,
      col sep=comma,
      x=rbsize,
      y=relerr,
    ]{pics/data/plot_rSVD-likeXs_povs10_qpow0.dat};

    \addplot[rSVDlikep30q0] table [
      header=true,
      col sep=comma,
      x=rbsize,
      y=relerr,
    ]{pics/data/plot_rSVD-likeXs_povs30_qpow0.dat};

    \addplot[rSVDlikep30q1] table [
      header=true,
      col sep=comma,
      x=rbsize,
      y=relerr,
    ]{pics/data/plot_rSVD-likeXs_povs30_qpow1.dat};
\end{groupplot}


\node[
    left = 0.9cm of err_plots c1r1.west,
    align=left,
    anchor=base,
    rotate=90,
  ] {basis gen. time};

\node[
  left = 0.9cm of err_plots c1r2.west,
  align=left,
  anchor=base,
  rotate=90,
] {rel.\ red.\ error};

\node[
    below = 1cm of err_plots c3r2.south east,
    align=left,
    draw=black,
    anchor=north east,
  ] {%
  \scriptsize
  \ref*{leg0}
  };
\end{tikzpicture}%
  \caption{Basis generation times and relative reduction errors (\ref{eqnrelerr}), rSVD-like}%
  \label{fig:methcomp1}%
\end{figure}

\vspace{-3mm}
\label{sec:numerics}
We apply the developed randomized, structure-preserving methods to a 2D linear wave equation model. The initial boundary value problem for the unknown $u(t, \fxi)$ with the spatial variable $\fxi := (\xi_1, \xi_2) \in \varOmega := (0, 1) \times (0, 0.2)$ and the temporal variable $ t \in \It(\fmu) := [\tInit, \tEnd(\fmu)]$, reads
\begin{align*}
u_{tt}(t, \fxi) &= c^2 \Delta u(t, \fxi) &&\textrm{in } \Omega\times\It(\fmu)\\
u(0, \fxi) &= u^0(\fxi) := h(s(\fxi)), \quad u_t (0, \fxi) = v^0(\fxi) := 0 &&\textrm{in } \Omega,\\
u(t, \fxi) &= 0 &&\textrm{in } \partial\Omega \times\It(\fmu), 
\end{align*}
with \vspace{-3mm}
\begin{align*}
  &s(\fxi) = 10\cdot\left(\xi_1 - \frac{1}{2}\right),&
  &h(s) =
  \begin{cases}
  1 - \frac{3}{2} |s|^2 + \frac{3}{4} |s|^3, & 0 \leq |s| \leq 1 \\
  \frac{1}{4} (2 - |s|)^3, & 1 < |s| \leq 2  \\
  0, & |s| > 2.
  \end{cases}
\end{align*}
We choose $\tInit = 0,\, \tEnd(\fmu) = 2/\fmu$ and as parameter (vector) $\fmu = c \in[1,2]$. Spatial discretization via central finite differences leads to the Hamiltonian system
\begin{align}\label{wave_discr}
  \ddt \fx(t; \fmu)
= \JtN \grad[\fx] \Ham(\fx(t; \fmu); \fmu) =\JtN \fA(\fmu)\fx, \quad  \fx(\tInit; \fmu) = \fxInit(\fmu)
\end{align}
with \vspace{-3mm}
\begin{align*}
  &\fxInit(\fmu) = [u^0(\fxi_1)); ...;u^0(\fxi_N));\Z{N\times 1}  ],&
  &\fA(\fmu) = \begin{pmatrix}
\fmu^2(\fD_{{\xi_1}{\xi_1}}+\fD_{{\xi_2}{\xi_2}}) \ &\Z{N} \\
\Z{N}\  & \I{N}
\end{pmatrix},
\end{align*}
where
$\{ \fxi_i \}_{i=1}^N \subset \varOmega$ are the grid points and the positive definite matrices
$\fD_{{\xi_1},{\xi_1}},$ $\fD_{{\xi_2}{\xi_2}} \in \R^{N \times N}$ denote the three-point central difference approximations in $\xi_1$-direction and in $\xi_2$-direction.
The domain is discretized equidistantly with 1000 grid points in $\xi_1$-direction and 20 points in $\xi_2$-direction which results in $N=1000 \cdot 20 = 20000$ grid points in total. 
The corresponding Hamiltonian reads $\Ham(\fx; \fmu) = \frac{1}{2}\rT\fx\fA(\fmu)\fx $. Temporal discretization is achieved with the implicit midpoint rule  and $n_t = 1000 $ equidistant time steps.  This results in different time step sizes for different parameters. In \Cref{fig:methcomp1,fig:methcomp2}, we present basis generation times and the relative reduction error 
\vspace{-1mm}
\begin{equation}\label{eqnrelerr}
\fe_\text{rel}(\fmu) = \sqrt{\sum\limits_{i = 0}^{n_t}||\fx_i(\fmu)  - \fV\fxri(\fmu)||_2^2}\Bigg{/}\sqrt{\sum\limits_{i = 0}^{n_t}||\fx_i(\fmu)||_2^2}, 
\end{equation} 
with $  \fx_i(\fmu), \fxri(\fmu), i = 0,.., n_t$ the iterates of the full model time-stepping and reduced model time-stepping, in dependence on the basis size. All results are averaged over 10 random parameters $\fmu \in \mathcal{P}.$ The results of the randomized methods are additionaly averaged over 5 runs for different random sketching matrices $\fOmega$. The snapshot matrix consisting of $\ns = 2000$ snapshots is computed from the parameters $\fmu_1 = 1, \fmu_2 = 2.$ The rcSVD is compared with three different versions of the cSVD: For 'cSVD full' a full SVD of $\fXsc$ is computed and then truncated, for 'cSVD with svds' the Matlab function svds is used that computes only the first $k$ singular vectors and values, for 'cSVDev' the first $k$ eigenvectors $\fv_1,..., \fv_k$ and eigenvalues $\lambda_1,...,\lambda_k$ of $(\fXsc\rT)\fXsc$ are computed and $\fU_\textrm{C} =\fXsc\ [\fv_1/\sqrt{\lambda_1},..., \fv_k/\sqrt{\lambda_k}]$ is set. We observe that the basis generation times are strongly reduced by randomization. 
With one power iteration and some oversampling (almost) the same reduction error is obtained for the rcSVD and the rSVD-like\texttt{KXs}, compared to its classical version. With the other two versions of the rSVD-like either a competitive error for only small basis sizes (rSVD-like\texttt{K}) or only large basis sizes (rSVD-like\texttt{Xs}) is obtained. Overall, the randomized techniques are able to reduce the computational time for generating a basis yielding a certain error.
For example an accuracy of about $10^{-4}$ is achieved with about 80 basis vectors by both, the cSVD and the rcSVD with $\qpow = 1$ and $\povs =10$. But the basis generation of the rcSVD is about 2.5 times faster. One power iteration significantly improves the reduction error, but leads to higher runtimes. Oversampling with $\povs=10$ only slightly effects error and runtime. 
\vspace{-7mm}
\vspace{-15mm}
\vspace{-7mm}
With about 100 basis vectors, an accuracy of $10^{-4}$ is reached for both, the SVD-like and the rSVD-like\texttt{KXs} with one power iteration and $\povs = 30$. But the computational costs are more than 3 times less with the randomized approach. Again, one power iteration significantly improves the reduction error, but also increases the runtimes.  
\vspace{-4mm}
\section{Conclusion and Outlook}\vspace{-3mm}
\label{sec:conclusion}
In our work, we have shown that randomized matrix factorizations can be used for the structure-preserving basis generation for symplectic MOR.  
The newly presented methods, the rcSVD and the rSVD-like decomposition, both provide very accurate approximations and lead to significant computational speed-ups compared to their classical versions. 
Future work will deal with the question how a randomization of the real canonical form can be obtained for randomizing the SVD-like decomposition instead of using a randomized Schur decomposition. 
\vspace{-3mm}
\subsubsection{Acknowledgements}
Supported by Deutsche Forschungsgemeinschaft (DFG, German Research Foundation) Project No. 314733389, and under Germany's Excellence Strategy -
EXC 2075 – 390740016. We acknowledge the support by the Stuttgart Center for Simulation Science (SimTech).
\vspace{-4mm}

\bibliographystyle{splncs04}
\bibliography{references_rand_sympl_mor}

\end{document}